\makeatletter
\expandafter\def\csname ver@ifluatex.sty\endcsname{}
\makeatother

\documentclass[12pt]{amsart}

\usepackage{amsmath,amssymb,amsfonts,amsthm,amscd}
\usepackage{mathtools}
\usepackage[mathscr]{eucal}
\usepackage{latexsym,graphics,enumerate,setspace,color,nicefrac,xfrac}
\usepackage[left=1in,right=1in,top=1in,bottom=1in]{geometry}
\usepackage{tikz}
\usepackage[colorlinks=true,linkcolor=red,urlcolor=blue,citecolor=gray]{hyperref}

\newcommand{\RR}{\mathbb{R}}

\newcommand{\be}{\begin{equation}} 
\newcommand{\bea}{\begin{eqnarray}}
\newcommand{\ee}{\end{equation}} 
\newcommand{\eea}{\end{eqnarray}}
\newcommand{\ba}{\begin{aligned}}
\newcommand{\ea}{\end{aligned}}

\newtheorem{corollary}{Corollary}[section]
\newtheorem{definition}{Definition}[section]

\theoremstyle{remark}
\newtheorem{remark}{Remark}[section]
\newtheorem{example}{Example}[section]

\renewcommand{\geq}{\geqslant}

\renewcommand{\leq}{\leqslant}

\def\d{{\textnormal{d}}}

\def\bcc{{\mathbf c}}
\def\bu{{\mathbf u}}
\def\bubar{\overline{\bu}}
\def\ubar{\overline{u}}

\def\bv{{\mathbf v}}
\def\sfv{{\mathsf v}}
\def\sfvp{{\mathsf v}'}
\def\rhop{\rho'}
\def\rhov{\rho\hspace*{0.02cm}\sfv}
\def\bx{{\mathbf x}}
\def\bxp{\bx'}
\def\bq{{\mathbf q}}

\def\bphi{\boldsymbol{\phi}}
\def\buE{{\bu}^{\scalebox{0.8}{\text{\tiny Euler}}}}
\def\bffE{{\bff}^{\scalebox{0.8}{\text{\tiny Euler}}}_j}

\def\buplus{\bu{\scriptscriptstyle +}}
\def\buminus{\bu{\scriptscriptstyle -}}
\def\xplus{x{\scriptscriptstyle +}}
\def\xminus{x{\scriptscriptstyle -}}
\def\bw{{\mathbf w}}
\def\bff{{\mathbf f}}

\def\scrC{{\mathscr C}}

\def\eps{\epsilon}

\def\seteta{{\mathscr E}} 
\def\setU{{\mathscr W}} 
\def\setUeps{\setU} 

\def\bxp{{\bx}'}

\def\dx{\d\bx}
\def\dxp{\d\bxp}

\def\dt{{\d}t}
\def\ds{{\d}s}
\def\etac{\eta_\bcc}
\def\myI{{\mathcal I}}
\def\Ieta{\myI_\eta}
\def\Ietac{\myI_{\eta_{{}_\bcc}}}
\def\Ietau{\myI_\eta(\bu)}
\def\Ietaubar{\myI_\eta(\bubar)}

\def\divbff{\sum_{j=1}^d \bff_j(\bu)_{x_j}}
\def\divbffbar{\sum_{j=1}^d \bff_j(\bubar)_{x_j}}

\numberwithin{equation}{section}

\begin{document}

\title[Variational formulation of hyperbolic conservation laws]{Variational formulation\\of hyperbolic conservation laws}

\author{Eitan Tadmor}
\address{Department of Mathematics and IPST, University of Maryland, College Park.}
\email{tadmor@umd.edu}

\date{\today}



\thanks{Research was supported  by  ONR grant N00014-2412659 and NSF grant DMS-2508407. I am also grateful for the hospitality of the Laboratoire Jacques-Louis Lions (LJLL) at Sorbonne University, where part of this work was completed}

\dedicatory{In memory of Peter Lax}

\begin{abstract}
Entropy functions played a key role in the development of mathematical theory for hyperbolic conservation laws. 
The notion of entropy, which is intimately connected with symmetry, is an extension \emph{imposed} on nonlinear systems of conservation laws. In this context, Friedrichs raised the question whether the assumed symmetries can also be derived. 
We introduce a variational formulation that addresses Friedrichs' question: an entropy function is derived as an extremal object, from which we deduce, rather than impose, the maximum entropy production principle of Dafermos.
\end{abstract}

\maketitle

\qquad  ``There is no theory for the initial value problem for compressible flows
in two 

\qquad   space dimensions once shocks show up, much less in three space dimensions. 

\qquad This is a scientific scandal and a challenge.'' \hspace*{2.1cm} P. D. Lax \cite{Lax2007}

\setcounter{tocdepth}{1}
\tableofcontents

\addtocontents{toc}{\protect\setcounter{tocdepth}{0}}
\section*{Dedication} 

Peter Lax was one of the preeminent mathematicians of the second half of the twentieth century and the leading  ambassador for modern applied mathematics \cite{Lax1989}. He served as a  role model for the generations of mathematicians who followed him, myself included.\newline
To grasp the scale of his mathematical legacy, one can recall the  concepts that bear his signature: the Lax equivalence theorem, the Lax-Milgram lemma, the Hopf-Lax formula, the Lax-Friedrichs and Lax-Wendroff methods, the Glimm-Lax memoir, Lax-Levermore theory, Lax pairs, Lax-Phillips scattering theory, the Lax entropic shock conditions,  the HLL Riemann solver, not to mention his pioneering works on the Riemann problem in hyperbolic conservation laws, and  on Fourier integral operators. A deeper look  into Lax's contributions can be found in \cite{Lax2005}.

I first met Peter during a trip to Northern California in the summer of 1981. I was a young postdoc from Caltech, giving a seminar at Stanford. When I mentioned that I would be visiting New York later that summer, he offered, entirely on the spot: ``Why don't you stay in my New York apartment?'' That ``you''  included our daughter, who was three years old at the time.
It remains one of those unforgettable life events that exemplified Peter's generosity. When I later shared this memory with his son, Dr. Jim Lax, he replied, ``My parents were like that. My friends referred to their apartment as Hotel Lax.''\newline
This was quintessentially Lax. He was not merely an iconic mathematical figure for my generation; he was an extraordinary human being. He was an admirer of John von Neumann \cite{Lax2014} and, as he himself noted, not the first Hungarian-born to serve as President of the American Mathematical Society. 
Both were part of the Hungarian group of scientists affectionately  referred to as ``The Martians" by their peers \cite{wiki:Martians}.\newline 
Peter had style and class, both in his life and in his mathematics. While his impact cannot be adequately captured in a few words, his legacy will undoubtedly endure.

\addtocontents{toc}{\protect\setcounter{tocdepth}{1}}
\section{Introduction --- a question of Friedrichs about symmetry}

In 1979, Friedrichs concluded his von Neumann lecture with the following question \cite{Fri1980}:

\smallskip
 ``\emph{For the systems of equations I have discussed here, the symmetry feature was a derived property. Now, in many branches of physics $\ldots$ symmetries play a fundamental role, but all these
symmetries --- as it seems to me --- are assumed and not derived. I now wonder whether or
not $\ldots$ symmetries can also be derived
from the overdeterminancy of basic conservation equations?}''

\smallskip
This paper describes an attempt to answer Friedrichs' question in the context of hyperbolic systems of conservation laws.\footnote{Dafermos  concludes his sketch of the early history of hyperbolic conservation laws, writing \cite[Introduction]{Daf2016},
``The next major milestone $\ldots$
is the landmark paper by Lax \cite{Lax1957}, which coins the term ``hyperbolic conservation law'' and launches the field as a new principal branch in the theory of partial differential equations. This was accomplished by distilling, generalizing and formalizing the raw material that had accumulated over the years  $\ldots$  It is fair to say that Lax's paper set the direction for the development of the field of hyperbolic conservation laws over the past fifty years.''
}
\begin{equation}\label{eq:CLs}
\bu_t + \sum_{j=1}^d \bff_j(\bu)_{x_j}=0, \qquad \bu(t,\cdot):\RR_\bx^d \mapsto \RR^N, \quad t\geq 0, \  \bx=(x_1,\ldots, x_d) \in \RR^d.
\end{equation}
 Here, $\bu=(u^1,\ldots, u^N)^\top$ is an $N$-vector of conserved quantities  which are balanced by flux vectors, $\bff_j(\bu)=(f_j^1(\bu),\ldots,f_j^N(\bu))^\top: \RR^N \mapsto \RR^N$. Integrating \eqref{eq:CLs} yields the conservative weak form for the $N$ components of the equations, $u^k_t + \text{div}\,\bff^k(\bu)=0$,
 \begin{equation}\label{eq:conservative}
 \int \limits_\Omega u^k(t,\cdot)\dx\Big|^{t=t_2}_{t=t_1}+\int \limits_{t=t_1}^{t_2}
  \int \limits_{\partial\Omega}\bff^k\cdot {\mathbf n}\, \d S=0, \quad 
 \bff^k(\bu)=\big(f^k_1(\bu),\ldots, f^k_d(\bu)\big), \ \ k=1,\ldots,N.
 \end{equation}
It states that the change over  time in the  amount of matter inside an arbitrary spatial domain $\Omega \subset \RR^d$ is solely due to the flux of that matter across the boundary of $\Omega$. We recall that \eqref{eq:conservative} is equivalent to the usual notion of weak solutions, where \eqref{eq:CLs} is interpreted in the distributional sense, \cite{Lax1973}.\newline
 Differentiation of \eqref{eq:CLs} puts it in the form of a quasi-linear system
 \begin{equation}\label{eq:quasi-linear}
 \bu_t +\sum_{j=1}^d A_j(\bu)\bu_{x_j}=0, \qquad A_j(\bu):=\frac{\partial \bff_j(\bu)}{\partial \bu}.
 \end{equation}
  The symmetry  mentioned in Friedrichs' question  is connected with the notion of an entropy. A nonlinear scalar function, $\eta=\eta(\bu)$, is an entropy function associated with \eqref{eq:CLs}, if there exists an entropy flux, $\bq^\eta(\bu)=(q_1(\bu),\ldots,q_d(\bu))^\top$, such that\footnote{$\nabla\zeta(\bu)$ is the \emph{column} gradient vector $(\zeta_{u^1}(\bu),\ldots, \zeta_{u^N}(\bu))$ and $\nabla\zeta(\bu)^\top$ is the \emph{row} gradient vector.}
  \begin{equation}\label{eq:entropy-pair}
  \nabla \eta(\bu)^\top A_j(\bu)=\nabla q_j(\bu)^\top, \qquad j=1,\ldots,d.
  \end{equation}
It follows that the conservation law \eqref{eq:CLs} admits an extension in terms of the entropy pair $(\eta,\bq^\eta)$:
 \[
 \Big\langle \nabla\eta(\bu),\bu_t+\sum_{j=1}^d \bff_j(\bu)_{x_j}\Big\rangle = 
 \eta(\bu)_t + \sum_j \nabla \eta(\bu)^\top A_j(\bu)\bu_{x_j} = \eta(\bu)_t + \nabla_\bx\cdot \bq^\eta(\bu).
 \]
 Of course, such an extension holds for any linear combination of the components of $\bu$; the essence, therefore, is  imposing the compatibility relation \eqref{eq:entropy-pair} for  nonlinear entropy functions. Let $\scrC(\RR^N,\RR)$ denote the cone of strictly convex functions, then one is led to the requirement that for all $\eta\in \scrC$
 \begin{equation}\label{eq:entropy-ineq}
 \eta(\bu)_t + \nabla_\bx\cdot\bq^\eta(\bu)\leq 0.
 \end{equation}
 This is the celebrated \emph{entropy inequality} imposed as a selection criterion of physically relevant weak solutions for \eqref{eq:CLs} \cite{Lax1971,Kru1970}.
  
  In \cite{FL1971}, Friedrichs and Lax 
 observed that the entropy compatibility requirement \eqref{eq:entropy-pair} is equivalent to symmetry.  This follows from differentiation of \eqref{eq:entropy-pair}  (here and below $D^2\zeta(\bu)$ denotes the Hessian $\displaystyle (D^2\zeta)_{\alpha,\beta}=\frac{\partial^2\zeta}{\partial u^\alpha\partial u^\beta}$),
 \[
  (D^2\eta) A_j + T_j = D^2q_j, \qquad \big(T_j\big)_{\alpha\beta}:=\Big\langle \nabla \eta,\frac{\partial^2 \bff_j}{\partial u^\alpha\partial u^\beta}\Big\rangle.
 \] 
 Hence, since the Hessians $D^2q_j$ and the tensors $T_j$ are symmetric, the existence of an entropy implies that its Hessian symmetrizes the Jacobians 
 \begin{equation}\label{eq:symmetry}
  D^2\eta(\bu)A_j(\bu) = \big(A_j(\bu) D^2\eta(\bu)\big)^\top \equiv A^\top_j(\bu)D^2\eta(\bu), \quad j=1,\ldots,d.
 \end{equation}
 Conversely, if the conservative system \eqref{eq:entropy-pair} is symmetrizable by a positive Hessian,  $(D^2\eta) A_j=A^\top_j (D^2\eta)$, then there exists entropy flux $\bq^\eta=(q_1,\ldots, q_d)^\top$ such that \eqref{eq:entropy-pair} holds and $(\eta,\bq^\eta)$ forms a (convex) entropy pair.
 
 \smallskip
We conclude  this Introduction with two comments which elaborate on the interplay between symmetry and entropy.
\subsection{Imposing the existence of entropy symmetrizer $N\geq 3$}  The symmetry condition \eqref{eq:symmetry}  clearly holds in the scalar case: When $N=1$,  every convex $\eta(u)$ serves as a convex entropy. In  particular, 
Kru\v{z}kov's one-parameter family of   entropy functions depending  on a ``dual'' scalar variable $c$ (or equivalently, the kinetic velocity variable in the kinetic formulation of \cite[Corollary 1]{LPT1994}), 
\begin{equation}\label{eq:Kruzkov}
\eta_c(u)=|u-c|, \qquad c\in \RR,
\end{equation}
 is  the main tool for development of existence, uniqueness and $L^1$-stability of scalar conservation laws, \cite{Kru1970} (see also the earlier work  \cite[\S16]{Vol1967}). In the case of $N=2$ equations, the symmetry \eqref{eq:symmetry} amounts to a second-order linear equation for $D^2\eta=D^2\eta(u^1,u^2)$. In this context we mention the family of Lax entropy pairs \cite{Lax1971}, and  special classes of $2\times 2$ systems with entropy   functions that arise as solutions for Euler-Poisson-Darboux equation, \cite{DiP1983},\cite[\S1.1]{LPT1994},\cite{CL2000}, \cite[\S3]{Che2005}, or  Goursat problem, \cite[\S9.3]{Ser1999}, \cite[Chap. XII]{Daf2016}. For $N\geq3$, the symmetry condition forms an over-determined system for the  $N(N-1)/2$ entries of $D^2\eta$.   
 The  $3\times 3$ system 
  \begin{equation}\label{eq:Rod}
 \bu_t + \begin{bmatrix}u^2 & & \\ & u^3 & \\ & & u^1\end{bmatrix}\bu_x=0,
 \end{equation}
is an example for the class of ``completely non-conservative'' systems studied by   Rozhdestvenskii, \cite{Roz1959} (translated in \cite[\S7]{Roz1960}).
It is in this sense that one needs to impose the entropy symmetrizer condition \eqref{eq:symmetry}, at least for systems with $N\geq3$ equations.  

\subsection{Symmetry and conservation}The entropy Hessian in \eqref{eq:symmetry}, $A_0^S(\bu):=D^2\eta(\bu)>0$, puts the system \eqref{eq:quasi-linear} into Friedrichs' symmetric form, \cite{Fri1958},
 \begin{equation}\label{eq:symmetrized-system}
 A_0^S(\bu)\bu_t + \sum_{j=1}^d A_j^S(\bu)\bu_{x_j}=0, \qquad A_j^S(\bu):=A_0^S(\bu)A_j(\bu), \ \ j=1,\ldots,d.
 \end{equation}
The system is symmetric in the sense that $A_j^S, j=0,1, \ldots,d$ are symmetric, and  therefore hyperbolic in the sense that 
$\sum_{j=0}^d A_j^S \omega_j$ are diagonalizable with real  eigenvalues. However, symmetrization  comes at the expense of conservation: the symmetric  $A_j^S(\bu)$ need not be  Jacobian matrices which would enable to identify \eqref{eq:symmetrized-system} as a system of conservation laws.

In his seminal  1961 paper entitled ``An interesting class of quasilinear systems'' \cite{God1961}, Godunov identified a class of quasilinear equations 
that can be written in a form that is both symmetric and in conservation form:
\begin{equation}\label{eq:interesting}
\big(\nabla L_0(\bv)\big)_t + \sum_{j=1}^d \big(\nabla L_j(\bv)\big)_{x_j}=0, 
\end{equation} 
with a quasi-linear symmetric form expressed in terms of the corresponding Hessians
\begin{equation}\label{eq:interesting-symmetric}
(D^2_\bv L_0)\bv_t + \sum_{j=1}^d (D^2_\bv L_j)\bv_{x_j}=0.
\end{equation}
Godunov showed that the equations of compressible gas dynamics as well 
as other hyperbolic systems in mathematical physics admit the conservative form \eqref{eq:interesting} for a proper choice of $L_j=L_j(\bv)$ and variables $\bv$ that were worked out in each case.\newline
A final and decisive step, bridging the work of Godunov with that of Friedrichs and Lax, was taken in  1980 by Mock, \cite{Moc1980}. Given a convex  entropy function $\eta(\bu)$, one defines the \emph{entropy variables}, $\bv:=\nabla \eta(\bu)$; by convexity one can, at least locally, consider the inverse mapping which we express as  $\bu=\bu(\bv)$. Expressed in terms of these entropy variables, system \eqref{eq:CLs} keeps its conservative form
\begin{equation}\label{eq:entropy-variable-form}
\bu(\bv)_t + \sum_{j=1}^d \bff_j(\bu(\bv))_{x_j}=0.
\end{equation}
Moreover, the $\bv$-dependent fluxes are now perfect gradients of an entropy potential, $\psi_0(\bv)$, and the corresponding flux potentials, $\psi_j(\bv)$,
\[
\begin{split}
\bu(\bv)&=\nabla \psi_0(\bv), \qquad \psi_0(\bv):= \langle \bv,\bu(\bv)\rangle -\eta(\bu(\bv))\\
\bff_j(\bu(\bv))& = \nabla\psi_j(\bv), \qquad \psi_j(\bv):=\langle \bv,\bff_j(\bu(\bv))\rangle-q_j(\bu(\bv)).
\end{split}
\]
This yields the quasi-linear  form of \eqref{eq:entropy-variable-form} expressed in terms of the symmetric Hessians
\begin{equation}\label{eq:entropy-variable-quasi}
\big(D^2_\bv\psi_0\big) \bv_t + \sum_{j=1}^d \big(D^2_\bv\psi_j\big)\bv_{x_j}=0.
\end{equation}
Compared with the Friedrichs-Lax symmetrization in \eqref{eq:symmetry}, we observe that the entropy variables formulation in \eqref{eq:entropy-variable-form} symmetrizes the Jacobians $A_j$ ``on the right'', since \eqref{eq:symmetrized-system} yields
\begin{equation}\label{eq:symmetry-in-v}
 \frac{\partial \bff_j(\bu(\bv))}{\partial \bv} =A_j(\bu)(D^2\eta)^{-1}=\big(A^S_0\big)^{-1}A_j^S(\bu)\big(A^S_0\big)^{-1}.
 \end{equation}
  This type of ``symmetrization of the right'' using the entropy variables was the main tool in the construction of entropy conservative/stable schemes in \cite{Tad1987a, Tad1987b,Tad2003}.
 
\smallskip
Observe that the potentials that dictate the entropy variable formulation \eqref{eq:entropy-variable-quasi},  \mbox{$\{\psi_j(\bv), j=0,\ldots,d\}$}, coincide with the 
functionals,  \mbox{$\{L_j(\bv), j=0,\ldots, d\}$}, in \eqref{eq:interesting}. Thus, we conclude that  \emph{Godunov's  ``interesting class of quasi-linear systems''  encompasses \emph{all} hyperbolic conservation laws endowed with (at least one) convex entropy extension.}

\begin{example}[{\bf The Euler system for compressible gas dynamics}]
 When $N\geq 3$,  the existence of an entropy function is the exception rather than the rule. Nevertheless, many physically relevant systems are endowed with entropy functions.
 The primary example --- the one that in fact motivated much of the theoretical development over the years --- is the Euler system for compressible gas dynamics. It consists of  $N=d+2$ equations 
for the conservative vector function, $\bu=\buE$, which consists of density $\rho>0$, momentum $\rhov\in \RR^{d}$, and total energy $\displaystyle E:=\frac{1}{2}\rho|\sfv|^2 + \rho e$,
 \begin{equation}\label{eq:Euler-eq}
 \buE=\begin{bmatrix}\rho\\ \rhov\\ E\end{bmatrix}, \quad \bffE(\bu)=\begin{bmatrix}\rho\sfv_j\\ \displaystyle \rho\sfv_j\sfv+ p\delta_{ij}\\[0.1cm] \displaystyle  \sfv_j(E+p)\end{bmatrix}, \ \  j=1,\ldots, d, \quad  p=(\gamma-1)\rho e.
 \end{equation}
 It is endowed with a family of convex entropy functions of the specific entropy $S:=\ln(p\rho^{-\gamma})$, \cite{Har1983}
 \begin{equation}\label{eq:Euler-entropies}
 \eta(\buE)= -\rho h(S), \quad h'-\gamma h''>0, \ h'>0.
 \end{equation}
In particular, therefore, it implies the well-known symmetrizability of  the system of compressible Euler equations.\footnote{The spatial Jacobians are similar to  --- and therefore are symmetrizable into,  the particularly simple form 
 \[
 A_j(\bu)\sim \sfv_j{\mathbb I}_{N\times N} + 
 c\left[\begin{array}{c|c|c} 0 & {\mathbf e}_j^\top & 0 \\ \hline   {\mathbf e}_j & \displaystyle \vphantom {e^{e^{e^\top}}} {\Large 0}_{d\times d} \vphantom{\frac{a}{b}}& \\ \hline
 0 & & 0\end{array}\right], \qquad c:=\sqrt{\frac{\gamma p}{\rho}}.
 \]}
 The entropy variables for \eqref{eq:Euler-eq} corresponding to  $\displaystyle h(S)=\frac{\gamma+1}{\gamma-1}e^{S/(\gamma+1)}$ take the particularly simple form \cite[\S2]{Tad2003}
\[
\bv=(p\rho)^{-\frac{\gamma}{\gamma+1}}(E,-\rhov,\rho)^\top.
\]
Euler fluxes, $\bffE(\bu)$, have the distinctive feature of being homogeneous of degree 1 and their homogeneity is preserved  in their symmetric entropy-variable formulation. This is outlined in Appendix \ref{sec:appa}.
\end{example}
Besides canonical examples of 
Euler equations \eqref{eq:Euler-eq} and other systems of mathematical physics, there are other classes of hyperbolic systems with non-empty set of  entropy functions.
\begin{example}[{\bf Symmetric systems}]
We consider the class of hyperbolic systems \eqref{eq:quasi-linear} with symmetric Jacobians,
\[
A_j(\bu):= \frac{\partial \bff_j(\bu)}{\partial \bu}=A^\top_j(\bu), \qquad  j=1,\ldots,d.
\]
These are precisely the systems with fluxes  $\bff_j(\bu)$ induced by potentials,  
\[
\zeta_j(\bu):=\int^\bu \bff_j(\bw)\cdot\d\bw,\qquad j=1,\ldots, d
\]
such that $\nabla \zeta_j(\bu)=\bff_j(\bu)$.
Then the quadratic ``energy'' $\eta(\bu)=\frac{1}{2}|\bu|^2$ is an entropy function  \cite{God1961} and hence 
$\eta_\bcc(\bu)=\frac{1}{2}|\bu-\bcc|^2$ is a family of entropy functions parameterized by $\bcc\in \RR^N$.
 It follows that in the particular case of one-dimensional  systems,  
$\zeta(\bu)$ symmetrizes one-dimensional symmetric systems,  \cite{Tad1987c}
\[
\langle \nabla\zeta(\bu), \bu_t + \bff(\bu)_x\rangle
=\zeta(\bu)_t + \Big(\frac{|\bff(\bu)|^2}{2}\Big)_x, \qquad
\zeta(\bu)=\int^\bu \bff(\bw)\cdot\d\bw
\]
Though $\zeta(\bu)$ need not be convex, one may consider
the one parameter family of convex entropy functions $\eta_\lambda(\bu)$
\[
\eta_\lambda(\bu)=\frac{1}{2}|\bu|^2- \lambda \zeta(\bu), \qquad  \lambda A(\bu)< {\mathbb I}. 
\]
\end{example}

\section{A variational principle}
 We established the relation between symmetry and entropy in the context of nonlinear conservation laws:  the existence of  an entropy function  amounts to the assumption of existence of a symmetrizer,  \eqref{eq:symmetry}, or ``symmetrization on the right'' in an entropy variables formulation \eqref{eq:entropy-variable-form}.  This assumption of symmetry/existence of an entropy is  \emph{imposed} on $N$-systems of conservation laws when $N\geq 3$.  We now reformulate the question  of Friedrichs as follows: 

\smallskip
\emph{Is it possible to derive the existence of entropy function(s) rather than  merely assuming} 

\emph{their existence?}

\smallskip\noindent
We aim to answer this question in terms of the following  variational principle. It involves state vector functions, $\bu(t,\cdot) \in \setU \subset \text{BV}\cap L^\infty$, and a set of strictly convex  functions acting on these states, $\eta \in \seteta \subset \scrC(\RR^N,\RR)$, which we refer to as  \emph{observables}. The precise nature of the states admitted into $\setU$, and the observables  admitted into $\seteta$, will be explored below. At this stage we postulate that 
 $\seteta$ should be an affine space so that each  $\eta\in \seteta$ generates the $N$-parameter family $\eta_\bcc \in \seteta$,
\[
\text{if} \ \ \eta\in \seteta \ \ \text{then} \ \ \etac(\bu):=\eta(\bu)+\langle \bcc, \bu\rangle \in \seteta \quad \text{for all} \ \ \bcc\in \RR^N.
\]
\begin{definition}[{\bf Variational solutions}]\label{def:var} We say that $(\eta, \bubar)$ is a \emph{variational solution} of \eqref{eq:CLs} if $\eta$ belongs to a non-empty set of observables $\seteta$, and $\bubar\in\setU$ is a minimizer of an  action functional, $\Ietac(\bu)$, ``observed'' by $\etac$, such that for arbitrary time interval $(t_1,t_2)\subset (0,\infty)$ and smoothly bounded domains 
$\Omega \subset \RR^d$,
\begin{equation}\label{eq:variational-principle}
\Ieta\big(\bubar; (t_1,t_2)\big) \leq \Ietac\big(\bu;(t_1,t_2)\big)\ \ \textnormal{for all} \   \left\{\begin{array}{ll} \bu(t,\cdot)\in \setU, & t\in (t_1,t_2)\\  \bu(t,\cdot)=\bubar(t,\cdot), & t\leq t_1
\end{array}\right\} \ \textnormal{and} \   \bcc \in \RR^N.
\end{equation}
The action functional, $\Ieta$, is given by
\begin{equation}\label{eq:action}
 \Ietau  =\Ieta\big(\bu; (t_1,t_2)\big):=\int \limits_{t=t_1}^{t_2} \int \limits_\Omega \Big\langle \nabla \eta(\bu), \bu_t + \divbff\Big\rangle_{\!\bphi} \,\dx\dt.
 \end{equation}
In particular, the corresponding action $\Ietac$ is defined by the same formula with $\eta$ replaced by $\etac$.
\end{definition}
  
\begin{remark}
The space  $BV\cap L^\infty$  is a natural candidate  for a suitable  class of  solutions  of \eqref{eq:CLs} which must admit the emergence of jump discontinuities, both for multi-dimensional equations \cite{Kru1970} and for one-dimensional systems \cite{Gli1965,Ser1999,Bre2000,Daf2016} (with proper restrictions, e.g.,   \cite{Jen2000}). Consequently, the integrand of $\Ietau$ involves  non-conservative products which are interpreted in the sense of \cite{DLM1995}.  
Recall that given a vector function ${\mathbf b}(\bu): \RR^N \mapsto \RR^N$, the non-conservative product  for $\text{BV}$ data $\bu$, 
denoted $[{\mathbf b}(\bu)^\top\bu_x]_{\bphi}$, is  a uniquely defined Borel measure which depends on a vector-valued Lipschitz path $\bphi:[0,1]\times \RR^N\times \RR^N$, so that $\bphi(s)=\bphi(s; \buminus,\buplus)$ is connecting possible jump discontinuities $\buminus=\bu(\xminus)$ to $\buplus=\bu(\xplus)$, while $s\in [0,1]$. Specifically, the non-conservative product  at a regular point of discontinuity $x$ has Dirac mass with $\bphi$-dependent amplitude\footnote{This definition takes place in the $d=1$-dimensional case. The generalization for $\bx \in \RR^d$ in $d>1$ dimensions is implemented at regular points of jump discontinuities, so that their left and right limits, $\bu{\scriptscriptstyle\pm}=\bu(\bx{\scriptscriptstyle\pm})$, are well defined across surface of discontinuity; all other irregular points for $\text{BV}$ functions have $(d -1)$-dimensional measure $0$, \cite[\S9]{Vol1967}. This enables one to define $\big[{\mathbf b}(\bu)^\top \bu_{x_j}\big]_{\!\bphi}$, \cite[\S6.1]{DLM1995}.}
\[
\big[{\mathbf b}(\bu)^\top\bu_x\big]_{\!\bphi}(x) = \int_0^1 {\mathbf b}(\bphi(s))^\top\,\frac{\d \bphi(s)}{\ds} \ds.
\]
We find it convenient to use the notation
\[
\big\langle {\mathbf b}(\bu),\bu_x\big\rangle_{\!\bphi} = \big[{\mathbf b}(\bu)^\top\bu_x\big]_{\!\bphi}
\]
In particular, for given ${\mathbf b}:\RR^N \mapsto \RR^N$  we interpret  
\[
\big\langle {\mathbf b}(\bu),\bff_j(\bu)_{x_j}\big\rangle_{\!\bphi} =\big[\big({\mathbf b}(\bu)^\top A_j(\bu)\big)\bu_{x_j}\big]_{\!\bphi}.
\]
The DLM theory extends Volpert theory  \cite{Vol1967} for non-conservative products which is based on the straight path, $\bphi(s,\buminus,\buplus)=(1-s)\buminus+s\buplus$. 
Note that the variational statement takes  into account a  proper family of locally Lipschitz paths which  impact the value of the non-conservative products. 
\end{remark}

\subsection{Weak solutions}   A variational solution $\bubar$ should  minimize $\Ietac(\bu)$. In particular, since $\nabla\etac=\nabla\eta+\bcc$, we have
\[
 \Ietaubar \leq \Ietac(\bubar)=\Ietaubar
+\int \limits_{t=t_1}^{t_2} \int \limits_\Omega\Big\langle \bcc, \bubar_t + \divbffbar\Big\rangle_{\!\bphi} \,\dx\dt,
\]
hence
\[
\int \limits_{t=t_1}^{t_2} \int \limits_\Omega\Big\langle \bcc, \bubar_t + \divbffbar\Big\rangle_{\!\bphi} \,\dx\dt\geq 0 \ \ \text{for all} \ \ \bcc\in \RR^N.
\]
In the present case of a Jacobian coefficient matrix, the DLM product coincides with the distributional derivative of the corresponding flux,
$\bff_j(\bu)_{x_j}=\langle A_j(\bu),\bu_{x_j}\rangle_{\!\bphi}$,
 and is therefore  independent of the path $\bphi$.
Since the integrand is conservative, its integral is path independent, and we conclude
\[
\int \limits_{t=t_1}^{t_2} \int \limits_\Omega  \Big\{\ubar^k_t 
+ \text{div}\, \bff^k(\bubar) \Big\}\dx\dt=
 \int \limits_\Omega \ubar^k(t,\cdot)\dx\Big|^{t=t_2}_{t=t_1}+\int \limits_{t=t_1}^{t_2}
  \int \limits_{\partial\Omega}\bff^k(\bubar)\cdot {\mathbf n}\, \d S=
0, \ \ k=1,\ldots,N.
\]
\begin{corollary}\label{cor:var-is-weak} A variational solution is  a  weak solution of \eqref{eq:CLs}.
\end{corollary}
 The ``linearization'' argument of $\eta_\bcc$ for large $|\bcc|>\|\nabla\eta(\bu)\|_{L^\infty}$ which recovers $\bubar$ as a weak solution of the ``underlying'' conservation law, is  reminiscent of the recovery of scalar weak solutions, $u$, by  linearization of Kru\v{z}kov's entropies \eqref{eq:Kruzkov}, $\eta_c=|u-c|$   for large $c>\|u(t,\cdot)\|_{L^\infty}$. 
 
 \medskip
 It is clear that the closure of $\seteta$ under linear perturbations, $\big\{\etac  :   \bcc\in \RR^N\big\}$, is in fact equivalent to having  the class of variational solutions, $\bubar$, include weak solutions of \eqref{eq:CLs}.  In what follows, we shall need to consider perturbations of such weak solutions. Accordingly, we  let the set of admissible states,   $\setUeps$,  include  arbitrarily small $\eps$-perturbations of  weak solutions,
\begin{equation}\label{eq:setUeps}
\setUeps(t_1,t_2)=\left\{ \bu+\delta\bu \ : \
\bu \ \text{is a weak solution of \eqref{eq:conservative}};  \  \left\{\begin{array}{l}
 \vspace*{0.1cm}
 \displaystyle \delta\bu \in C_0^\infty\big((t_1,t_2), \Omega\big) \\
   \displaystyle \|\delta\bu\|_\infty< \eps \end{array}\right\}\right\}.
\end{equation}
The notion of a variational solution takes the following form.
 
 \medskip\noindent
 {\bf Variational solutions revisited. \!\!\!I}. $(\eta,\bubar)$ is a variational solution of \eqref{eq:CLs} if $\eta$ belongs to the non-empty set of observables $\seteta$,  and  $\bubar$ is a weak solution which is a minimizer of $\Ieta(\bu)$ observed by
$\eta$, such that for all perturbed $\bu \in \setUeps$ with arbitrarily small $\eps$, 
\begin{equation}\label{eq:variational-principle-weak}
\Ieta\big(\bubar; (t_1,t_2)\big) \leq \Ieta\big(\bu;(t_1,t_2)\big) \ \textnormal{for all} \ \   \left\{\begin{array}{ll} \bu(t,\cdot)\in \setUeps(t_1,t_2), & t\in (t_1,t_2)\\  \bu(t,\cdot)=\bubar(t,\cdot), & t\leq t_1.
\end{array}\right\}.
\end{equation}
The variational formulation \eqref{eq:variational-principle-weak} is based on a ``push forward'' action functional $\Ietau$. In particular, the admissible  states, $\bu(t,\cdot)$ --- and therefore the variational solution $\bubar(t,\cdot)$ --- are anchored in the initial data, $\bu(0,\cdot)=\bu_0$, prescribed at $t_1=0$. Thereafter, the states $\bu(t,\cdot), \ t>t_1>0$ are perturbations of  weak solutions which evolve from  the variational solution $\bubar(t_1,\cdot)$.
 
\smallskip\noindent
\subsection{Derivation of entropy functions} We  compute the formal first variation in $\Ietau$ using smooth perturbations, $\|\delta\bu\|_\infty< \eps$, compactly supported inside $(t_1,t_2)\times \Omega$,
\[
\begin{split}
\Ieta&(\bu+\delta\bu)-\Ietau\\
 \approx&
\int \limits_{t=t_1}^{t_2} \int \limits_\Omega\left\{\Big\langle D^2\eta(\bu)\delta\bu, \bu_t+ \sum_{j=1}^d A_j(\bu)\bu_{x_j}\Big\rangle_{\!\bphi} 
+\Big\langle\nabla\eta(\bu),(\delta\bu)_t + \sum_{j=1}^d \big(A_j(\bu)\delta\bu\big)_{x_j}\Big\rangle_{\!\bphi}\right\}\,\dx\dt.
\end{split}
\]
The approximation on the right takes into account first-order perturbations in $\delta\bu$, and ignores terms of order $\mathcal O(\eps^2)$, given that   that $\bu$ has bounded amplitude in $BV \cap L^\infty$.\newline
 The two terms on the right that involve time derivatives cancel each other
\[
\begin{split}
\int \limits_{t=t_1}^{t_2} \int \limits_\Omega&\left\{\big\langle D^2\eta(\bu)\delta\bu, \bu_t\big\rangle_{\!\bphi} 
+\big\langle\nabla\eta(\bu),(\delta\bu)_t \big\rangle_{\!\bphi}\right\}\,\dx\dt\\
 &= \int \limits_{t=t_1}^{t_2}  \int \limits_\Omega
\big\langle \nabla\eta(\bu),  \delta\bu\big\rangle_t \,\dx\dt=
\int \limits_\Omega\big\langle \nabla\eta(\bu),  \delta\bu\big\rangle\dx\Big|_{t=t_1}^{t=t_2}=
0.
\end{split}
\]
For the remaining spatial terms we have
 \begin{equation}\label{eq:by-parts}
 \begin{split}
&\sum_{j=1}^d \int \limits_{t=t_1}^{t_2} \int \limits_\Omega
 \left\{\Big\langle D^2\eta(\bu)\delta\bu,  A_j(\bu)\bu_{x_j}\Big\rangle_{\!\bphi} + \Big\langle\nabla\eta(\bu), (A_j(\bu)\delta\bu)_{x_j}\Big\rangle_{\!\bphi}\right\}\,\dx\dt\\
  & = \sum_{j=1}^d\int \limits_{t=t_1}^{t_2} \int \limits_\Omega\left\{\Big[ (\delta\bu)^\top  D^2\eta(\bu)A_j(\bu)\bu_{x_j}\Big]_{\!\bphi} - \Big[(\delta\bu)^\top A^\top_j(\bu)D^2\eta(\bu)\bu_{x_j}\Big]_{\!\bphi}\right\}\,\dx\dt.
 \end{split}
 \end{equation}
 This follows  from the following key feature of the DLM-Volpert theory: whenever a non-conservative product $\big[{\mathbf b}(\bu)^\top\bu_x\big]_{\!\bphi}$ forms a ``perfect derivative'', it coincides with the standard distributional derivative and is independent of the path $\bphi$, \cite[Proposition 1.5]{DLM1995},\cite[\S13.2]{Vol1967}; in particular,
\[
\big\langle {\mathbf b}(\bu)_x,\bu\big\rangle_{\!\bphi} + \big\langle {\mathbf b}(\bu),\bu_x\big\rangle_{\!\bphi}= \Big[ \Big(\bu^\top \frac{\partial{\mathbf b}(\bu)}{\partial \bu}+{\mathbf b}(\bu)^\top\Big) \bu_x\Big]_{\!\bphi}  = \big\langle{\mathbf b}(\bu),\bu\big\rangle_{\!x}, 
\]
 is independent of the path $\bphi$, which justifies  the ``integration by parts'' on the second  integrand in \eqref{eq:by-parts}.\newline 
Since we may choose arbitrary smooth, compactly supported perturbations $\delta\bu$, the requirement that $(\delta_\bu \Ieta)(\bubar)=0$ implies
\begin{equation}\label{eq:symmetrize-bubar}
\big[D^2\eta(\bubar)A_j(\bubar)\bubar_{x_j}\big]_{\!\bphi}=\big[A^\top_j(\bubar) D^2\eta(\bubar)\bubar_{x_j}\big]_{\!\bphi},
\end{equation}
and since this is intended to hold independently of the choice of path, we are led to
\[
D^2\eta(\bubar)A_j(\bubar)=A^\top_j(\bubar) D^2\eta(\bubar), \quad j=1,\ldots,d.
\]
This is precisely the symmetrizability condition \eqref{eq:symmetry} restricted to the variational solution $\bubar$. At the formal level, if the variational principle is required to hold for a sufficiently rich class of weak solutions, then the stationarity condition forces the pointwise symmetrization relation \eqref{eq:symmetry}.

In this manner the notion of an entropy function is deduced  from the variational principle \eqref{eq:variational-principle} or \eqref{eq:variational-principle-weak}: $\eta\in \seteta$ is an  admissible observable  if its Hessian   symmetrizes the Jacobians, \eqref{eq:symmetry}, that is, if $\eta$ is an entropy function.
This answers Friedrichs' question in the sense that the notion of entropy is derived as a stationarity condition for the variational functional $\Ieta$.

\begin{remark}[{\bf Entropic vs.\ non-entropic systems}]\label{rem:entropic-systems}The notion of variational solution does not treat the set of observable  as an ``extension'' which augments the hyperbolic system of conservation laws, but as an intrinsic part of the solution. Since observables coincide with entropy functions,  the set of observables $\seteta$ is restricted. Therefore, we make a distinction between two classes of systems of equations:\newline
 Entropic systems which admit non-empty set $\seteta$ --- this is the ``interesting'' class of conservation laws of Godunov
which admit (at least one) entropy function. This is the class of systems discussed in this paper. 
In case there is  more  than one entropy function, we can fix one preferred  $\eta \in \seteta$ and we refer to $\bubar$ as a variational solution observed by that $\eta$. At the same time, we can discuss a variational solution observed by all $\eta\in \seteta$.
Then there is the other class of   systems of conservation laws, at least when $N\geq 3$,   which admit no entropy extension, i.e.,  where $\seteta$ is empty, as in \eqref{eq:Rod}. This class of ``non-entropic'' systems do not admit variational solutions. The study of variational solutions for such systems requires a different setting.
\end{remark}

We close by noting the following duality. The closure of the set of observables  with respect affine perturbations, $\etac\in \seteta$, enforces variational solutions  to be weak solutions, $\bubar\in \setUeps$, while  letting the  set of admissible states $\setUeps$ include local perturbations  of weak solutions, $\bu+\delta\bu\in \setUeps$,  enforces the  observables to be entropy 
functions, $\eta\in  \seteta$. In this context, the class of perturbations  considered in \eqref{eq:setUeps}, which consists of all 
$\{\delta\bu\}$ in  the $L^\infty$ $\eps$-ball, is a ``rich'' class.
It should possible to consider a restricted  class of perturbations and consequently, to further restrict the set of perturbed weak solutions $\setUeps(t_1,t_2)$ in \eqref{eq:setUeps}. A key question of interest  is  \emph{to identify a minimal set of perturbed weak solutions which will enforce \eqref{eq:symmetrize-bubar}}.

\section{The entropy inequality}
Fix $\eta\in \seteta$.  Since its Hessian  symmetrizes the Jacobians, $A_j(\bu)$, it is a part of entropy-entropy flux pair, $(\eta,\bq^\eta)$, such that \eqref{eq:entropy-pair} holds, and the variational functional, $\Ietau$  now amounts to 
\begin{equation}\label{eq:variational-entropy}
\begin{split}
\Ieta\big(\bu; (t_1,t_2)\big) &= \int \limits_{t_1}^{t_2} \int \limits_\Omega 
\big\{\eta(\bu)_t + \nabla_\bx\cdot\bq^\eta(\bu) \big\} \dx\dt\\
& = \int \limits_\Omega \eta\big(\bu(t,\bx)\big)\dx\Big|^{t=t_2}_{t=t_1}+\int \limits_{t=t_1}^{t_2}
 \int \limits_{\partial\Omega}\bq^\eta\big(\bu(t,\bx)\big)\cdot {\mathbf n}\, \d S\, \dt.
 \end{split}
\end{equation}
The notion of variational solution in \eqref{eq:variational-principle-weak} now reads as follows.

 \medskip\noindent
 {\bf Variational solutions revisited. \!\!\!II}. Given an entropy function $\eta\in \seteta$, then $(\eta,\bubar)$ is a variational solution of \eqref{eq:CLs} if  $\bubar$ is a weak solution which is a minimizer of $\Ieta(\bu)$  among all weak solutions $\bu$ which agree with $\bubar$ at $t= t_1$, 
\begin{equation}\label{eq:variational-entropy-weak}
\begin{split}
\int \limits_\Omega  \eta\big(\bubar(t,\bx)\big)&\dx\Big|^{t=t_2}_{t=t_1}+\int \limits_{t=t_1}^{t_2}
 \int \limits_{\partial\Omega}\bq^\eta\big(\bubar(t,\bx)\big)\cdot {\mathbf n}\, \d S\, \dt\\
  & \leq \int \limits_\Omega \eta\big((\bu+\delta\bu)(t,\bx)\big)\dx\Big|^{t=t_2}_{t=t_1}+\int \limits_{t=t_1}^{t_2}
 \int \limits_{\partial\Omega}\bq^\eta\big((\bu+\delta\bu)(t,\bx)\big)\cdot {\mathbf n}\, \d S\, \dt.
 \end{split}
\end{equation}
Note that \eqref{eq:variational-entropy-weak} does \emph{not} recover the entropy inequality 
\eqref{eq:entropy-ineq}, which in its weak form reads
\begin{equation}\label{eq:weak-entropy}
\int \limits_\Omega  \eta\big(\bubar(t,\bx)\big)\dx\Big|^{t=t_2}_{t=t_1}+\int \limits_{t=t_1}^{t_2}
 \int \limits_{\partial\Omega}\bq^\eta\big(\bubar(t,\bx)\big)\cdot {\mathbf n}\, \d S\, \dt\leq 0.
\end{equation}
Indeed, the settings for  \emph{entropy solutions} satisfying \eqref{eq:weak-entropy}, and for variational solutions satisfying \eqref{eq:variational-entropy-weak}, are different. Entropy solutions are realized by vanishing viscosity limits, 
$\displaystyle \bu=\lim_{\eps \rightarrow 0} \bu^\eps$, where
  \[
  \eta(\bu^\eps)_t + \nabla_\bx\cdot\bq^\eta(\bu^\eps) \leq \eps\Delta_\bx \eta(\bu^\eps), \qquad \bu^\eps_t + \sum_{j=1}^d \bff_j(\bu^\eps)_{x_j}=\eps\Delta_\bx\bu^\eps, \qquad \eps>0.
  \]
  Viscosity perturbations impose an instantaneous decrease of entropy. A celebrated result of Lax, \cite{Lax1971,Lax1973}, shows that this  identifies physically relevant shock discontinuities,  encoded by nontrivial  defect measure  $\displaystyle \lim_{\eps\downarrow 0}\big\{-\eps|\nabla_\bx\eta(\bu^\eps(t,\cdot))|^2\big\}< 0$.
The notion of variational solutions replaces viscosity perturbations, $\bu \mapsto \bu^\eps$, by local smooth perturbations, $\bu \mapsto \bu+\delta\bu \in \setUeps$,  which, a priori, need not decrease the entropy production. See Example \ref{exm:alignment} below. 

\begin{remark}
It would be interesting to develop an alternative theory of variational solutions based on the class of Lipschitz
smooth perturbations, $\setUeps=\{\bu_{\!\bphi}\}$, where the non-conservative product $[{\mathbf b}(\bu)^\top \bu_x]_{\bphi}$ is replaced by the limit of 
$\displaystyle {\mathbf b}(\bu_{\!\bphi})^\top \frac{\d \bu_{\!\bphi}}{\d s}$, \cite[\S4]{DLM1995}.
\end{remark}

 The entropy inequality, $\eta(\bu)_t+\nabla_\bx\cdot\bq^\eta(\bu)\leq 0$,
was set  as a selection principle to identify a unique, physically relevant solution, \cite[Definition 2]{Kru1970}.
This selection principle   completely settled the question of uniqueness  in the scalar case,  being  endowed with the ``rich'' family of \emph{all} convex entropy functions, \cite{Kru1970}. But most systems are not endowed with  rich family of entropies, and those in mathematical physics are  identified with one preferred entropy, usually the one driven by thermodynamic considerations.   Therefore, one is interested to address the  uniqueness question in the setting ``observed'' by a single preferred entropy function. As noted in Remark \ref{rem:entropic-systems}, this is the preferred setting of variational  solutions which emphasizes a single entropy  function (or at least a  restricted set  $\seteta$ of entropies). We mention in this context two canonical results.
Panov \cite{Pan1994} proved the uniqueness of quadratic entropy solutions for one-dimensional scalar conservation laws, $u_t+f(u)_x=0$, with convex flux $f(u)$. See also \cite{LOW2004, KV2019}. For arbitrary systems of entropic conservation laws, the use of a single \emph{relative entropy} function, initiated in \cite{Daf1979}, secures uniqueness for within the class of strong solutions\footnote{As noted by Dafermos \cite[\S5]{Daf2016} 
``It is remarkable that a single entropy inequality, with convex entropy, manages
to weed out all but one solution of the initial value problem, so long as a classical solution exists.''}. However, a selection principle based on a single  entropy inequality    fails in the general setting of  systems of conservation laws in $N\geq2$ dimensions, see \S\ref{sec:entropy-uniqueness-fails} below.

\section{Maximum entropy production}\label{sec:MaxEnt} 
The possible lack of uniqueness in the entropy-based selection principle based on a preferred single  entropy inequality \eqref{eq:entropy-ineq}, motivated Dafermos as early as 1973, \cite{Daf1973}, to formulate an \emph{entropy rate admissibility criterion}. Given an strictly convex entropy $\eta$, it seeks a solution $\bubar$ which maximizes the entropy dissipation \emph{rate} among all other weak solutions
\begin{equation}\label{eq:Dafermos-MaxEnt}
\frac{\d_+}{\dt} \int  \eta(\bubar(t,\bx))\dx \leq 
\frac{\d_+}{\dt} \int  \eta(\bu(t,\bx))\dx\ : \ \left\{\begin{array}{l}
 \bu \ \text{is a weak solution}\\
 \bu(\tau,\cdot)=\bubar(\tau,\cdot), \ \ \tau\leq t.
 \end{array}\right.
\end{equation}
That is, there is no weak solution $\bu(\tau,\cdot)$ which coincides with $\bubar(\tau,\cdot)$ for $\tau<t$, and for which the \emph{reverse} of \eqref{eq:Dafermos-MaxEnt} holds.\newline
This can be interpreted in the general context of maximization
of the entropy production, the so-called MaxEnt principle, that  was developed in different contexts during the second half of the  20th century, with contributions (and sometimes controversies) from I. Prigogine (1947), J. M. Ziman (1956),
E. T. Jaynes (1957), H. Ziegler (1963), Dewar [2007] and many others listed in, e.g.,  \cite{MS2006}\footnote{As noted in \cite{GLC2020}, Dafermos' MaxEnt criterion relates to  Ziegler's principle for closed thermodynamic systems rather than  Prigogine's MinEnt principle principle which applies to open thermodynamic systems}.
Dafermos proved that, in the context of hyperbolic conservation laws, his MaxEnt criterion  singles out the unique entropy solution in scalar equations and 1D $p$-system \cite{Daf1973} and  selects the admissible solution of 1D Riemann problem, at least with small initial variation \cite{Sev1990,Daf2009},
and it was utilized in a host of different applications, e.g., \cite{HHW2020, Kle2023}.

\subsection{A local version of MaxEnt}\label{sec:localMax-Ent} We  appeal to the variational formulation \eqref{eq:variational-entropy-weak}.
Given that the observable $\eta$ is an entropy function, \eqref{eq:variational-entropy-weak} identifies   $\bubar$ as a variational solution which maximizes the local entropy production \emph{rate}, in the sense that $\bubar$ is a minimizer among all weak solutions, $\bu$, 
\begin{equation}\label{eq:Dafermos-MaxEnt-local-weak}
\begin{split}
\frac{\d_+}{\dt}\int \limits_\Omega &\eta(\bubar(t,\bx))\dx +
 \int \limits_{\partial\Omega}\bq^\eta(\bubar(t,\bx))\cdot {\mathbf n}\, \d S \\
  & \leq 
\frac{\d_+}{\dt} \int \limits_\Omega \eta(\bu(t,\bx))\dx+
 \int \limits_{\partial\Omega}\bq^\eta(\bu(t,\bx))\cdot {\mathbf n}\, \d S\ : \ \left\{\begin{array}{l}
 \bu \ \text{is a weak solution}\\
 \bu(\tau,\cdot)=\bubar(\tau,\cdot), \ \ \tau\leq t.
 \end{array}\right.
 \end{split}
\end{equation}
This is a \emph{local} version of Dafermos'  MaxEnt criterion \eqref{eq:Dafermos-MaxEnt}, localized to arbitrary $\Omega \subset \RR^d$. It  is  not assumed, but is deduced directly from the variational principle:  dividing  \eqref{eq:variational-entropy-weak} by $t_2-t_1$ and letting $t_2\rightarrow t_1=t$ while noting that $\delta\bu$ vanishes as well.

At this stage we can summarize the three main consequences of our variational principle outlined in Definition \ref{def:var}, namely ---
 variational solutions are weak solution;
stationary observable must be an entropy function; and
if $\eta$ is an entropy function, the variational principle implies local MaxEnt principle.

\section{In search of a uniqueness selection principle}

\subsection{The entropy inequality as a selection principle}\label{sec:entropy-uniqueness-fails}
There is a large body of recent work that addressed the (non-)uniqueness question of entropy solutions. The  canonical system  addressed in this large body of work is the isentropic Euler equations, which consists of $N=d+1$ equations for $(\rho, \rhov)\in \RR_+\times \RR^d$, 
 \begin{equation}\label{eq:Euler-isentropic}
 \bu=\begin{bmatrix}\rho\\ \rhov\end{bmatrix}, \quad \bff_j(\bu)=\begin{bmatrix}\rhov_j\\ \displaystyle \rhov_j \sfv + p\delta_{ij}\vspace*{0.1cm}\end{bmatrix}, \ \  j=1,\ldots, d, 
 \end{equation}
with pressure $p=p(\rho)$.  The $\gamma$-law, $p(\rho)=\rho^\gamma$, corresponds to Euler system \eqref{eq:Euler-eq} with  constant specific entropy $S$. Here, the energy 
$\displaystyle \eta(\rho,\rhov):=\frac{1}{2}\rho|\sfv|^2+ \rho \int^\rho \frac{p(r)}{r^2}\d r$ is serving as the entropy function.\newline
The relative entropy-based uniqueness within the class of strong  
 rarefaction-based solutions  was proved in \cite{CC2007,FK2015}.
Following the initial breakthrough of De Lellis and Sz\'{e}kelyhidi in 
\cite{DS2010} which used convex integration to construct non-unique ``wild'' solutions for bounded data, a series of results of non-uniqueness with improved 
regularity followed with regular density \cite{Chi2014},  compactly supported data \cite{AW2021}, 2D Lipschitz and smooth initial data \cite{CDK2015, CKMS2021}, and even non-unique entropy conservative H\"{o}lder solutions \cite{GK2022}. General non-uniqueness for shock-based solutions was proved in \cite{CK2014,CK2018,MK2018} and was extended to the full Euler equations in \cite{KKMM2020}.  
These works make clear that the selection principle based on a single entropy inequality to single out a unique weak solution from among solutions for the multidimensional isentropic equations \eqref{eq:Euler-isentropic}  does not suffice  for uniqueness.
\newline 
Lax commented on this situation in his Gibbs lecture \cite{Lax2007} ``\emph{Just because we cannot prove that compressible
flows with prescribed initial values exist doesn't mean that we cannot compute them.}''
Numerical evidence for non-uniqueness of entropic solutions  of Euler equations 
was further investigated  in \cite{FMT2016}.

\begin{example}[{\bf Increase of entropy}]\label{exm:alignment}
Like most of the theoretical aspects in this field, the entropy inequality \eqref{eq:entropy-ineq} was motivated by entropic solutions of the compressible  Euler equations \eqref{eq:Euler-eq}, which are realized as vanishing viscosity limits of Navier-Stokes equations (NSe) (here $\nabla_S$ denotes the symmetric gradient $\nabla_{\!{}_{\!\text{S}}} \cdot=\frac{1}{2}\big(\nabla\cdot +(\nabla \, \cdot)^\top\big)$)
\begin{equation}\label{eq:Navier-Stokes}
\buE_t+\sum_{j=1}^d \bffE(\bu)= \begin{bmatrix}0 \\ \text{div} \ \mathbb{T}
\\[0.1em] \nabla \cdot (\mathbb{T} \sfv+\kappa \nabla \theta)\end{bmatrix}, \quad
\mathbb{T} := 2\mu  \nabla_{\!{}_{\!\text{S}}} \sfv + \lambda (\nabla \cdot \sfv) \mathbb{I}.
\end{equation}

The Navier-Stokes equations (NSe) are \emph{the} canonical ``vanishing viscosity'' perturbation of  Euler equations \eqref{eq:Euler-eq}. They
 yield the entropy inequality for the convex entropy\footnote{The system of Euler equations is equipped with the family of entropy functions \eqref{eq:Euler-entropies}.  When  extended to the system of Navier-Stokes equations \eqref{eq:Navier-Stokes}, however,  only the ``physical'' entropy, $\eta(\buE) = - \rho S$, survives as the one which puts the additional viscosity and heat fluxes into a symmetric negative-definite form \cite{HFM1986}.} $\eta(\buE)=-\rho S$
uniformly as $\mu,\lambda,\kappa > 0$,
\[
(-\rho S)_t +\nabla_\bx\cdot (-\rhov S) < 0, \qquad S=\ln(p\rho^{-\gamma})
\]

However, the decrease of entropy is not shared by all well-posed perturbations of Euler equations. Consider the Euler alignment system, driven by a radial communication kernel $0< k \leq 1$,
\begin{equation}\label{eq:Euler-alignment}
\buE_t+\sum_{j=1}^d \bffE(\bu)=
\begin{bmatrix}0\\ \displaystyle \tau\int k(|\bx-\bx'|)\big(\sfvp-\sfv\big)\rho\rhop\dxp\\[1em] \displaystyle -\tau\int k(|\bx-\bxp|)\big(2E-\sfvp\cdot\sfv\big)\rho\rhop\dxp\end{bmatrix} \quad \sfvp=\sfv(t,\bxp).
\end{equation}
The system admits global weak solutions \cite{CCT2026}. These are in fact strong solutions for sub-critical initial data: it was shown in \cite{Tad2025} that Euler alignment, viewed as a perturbation of compressible Euler system \eqref{eq:Euler-eq},
yields the \emph{reverse} inequality uniformly in $\tau >0$,
\[
(-\rho S)_t+\nabla_\bx\cdot(-\rhov S)>0.
\]
This reflects the increasing order, or emergence in  the Euler alignment system
\eqref{eq:Euler-alignment}, and it demonstrates that imposing the ``direction'' of an entropy inequality \eqref{eq:entropy-ineq} is connected with the class of perturbations
applied to the  underlying systems of  conservation laws.
\end{example}

\subsection{The MaxEnt condition as a selection principle}
 Hsiao \cite{Hsi1980} found a surprisingly simple
counterexample of 1D Riemann problem for the compressible Euler
equations, where Dafermos' MaxEnt principle fails to single out a unique solution within the class of piecewise smooth solutions,  when $\gamma > 5/3$.
It turns out that  Dafermos’ criterion fails as a uniqueness selection principle for a certain class of initial data for the multidimensional isentropic Euler equations  \cite{CK2014,CDK2015}, that is, there exist wild solutions with better entropy dissipation rate when compared with the self-similar solution
of the 1D Riemann problem extended to two dimensions. Hence
the local maximal dissipation criterion rules out what we would expect to single out as the the unique solution in this case.\newline
Alternative approaches to Dafermos' MaxEnt principle were offered by the 
least action principle in
\cite{GGKS2025} as well as measure-valued solutions realized by low Mach limits or vanishing viscosity limits \cite{Gal2023, GW2021}, but they fail to secure uniqueness.\newline 
The MaxEnt principle of Dafermos was also examined in the context of the larger class of   \emph{measure valued}  solutions \cite{DiP1985} for \eqref{eq:CLs},
where various variational-based selection criteria were examined to single out  a unique measure-valued solution, primarily for the isentropic Euler equations \eqref{eq:Euler-isentropic}.
 In \cite{FJL2025} it was shown that a local MaxEnt selection principle identifies measured valued solutions of \eqref{eq:Euler-isentropic}, thus proving a (local version of)
DiPerna's conjuncture \cite{DiP1985} for 
an alternative two-step  selection criterion was studied in \cite{FL2025},  improving  the earlier multi-step selection criterion in \cite{BFH2020}. A different criterion, the MaxVar criterion for measure-valued solutions, was studied in \cite{KMW2026}; this remains a work in progress.

\subsection{A variational formulation as a selection principle} 
Markfelder \cite{Mar2024} proved that the local MaxEnt \eqref{eq:Dafermos-MaxEnt-local-weak} fails for the isentropic equations \eqref{eq:Euler-isentropic} with $\gamma=2$.
The variational formulation \eqref{eq:variational-entropy-weak} is in fact more restrictive than the local MaxEnt criterion \eqref{eq:Dafermos-MaxEnt-local-weak}, in the sense that it requires minimization over local perturbations of weak solutions. Its global version reads
\[
\int  \eta\big(\bubar(t,\bx)\big)\dx\Big|^{t=t_2}_{t=t_1}
   \leq \int\eta\big((\bu+\delta\bu)(t,\bx)\big)\dx\Big|^{t=t_2}_{t=t_1} \quad \textnormal{for all} \ \ \left\{\begin{array}{l}
   \bu+\delta\bu \in \setUeps(t_1,t_2)\\
   \bu(t,\cdot)=\bubar(t,\cdot).\end{array}\right\}.
\]
It remains an open question to identify a suitable space $\setUeps$ which would provide a fully rigorous treatment for  the questions of a single entropy-based selection  criterion for existence and uniqueness  of variational solutions.

\appendix
\section{Conservation laws with homogeneous fluxes}\label{sec:appa} The compressible  Euler equations \eqref{eq:Euler-eq} have the distinctive property of  spatial fluxes which are homogeneous, $\bffE(\lambda\bu)=\lambda^\beta\bffE(\bu)$ of degree $\beta=1$.
The purpose of this short appendix is to note that 
 the homogeneous structure of Euler equations can be put in a \emph{symmetric form}, \cite{Tad1984}. We think that it could be further utilized  in the variational formulation of Euler equations.\newline 
 To this end we use the class of convex entropy functions (corresponding to $\displaystyle h(S)=e^{S/(\alpha+\gamma)}$)
\begin{equation}\label{eq:Euler-homo}
\eta(\buE)=-(p\rho^\alpha)^{\frac{1}{\alpha+\gamma}}, \qquad \alpha \geq 0
\end{equation}
which are homogeneous of degree $\displaystyle \beta_\alpha=\frac{\alpha+1}{\alpha+\gamma}$. Expressed in terms of the entropy variables, $\bv=\nabla\eta(\buE)$, the Euler equations  admit the symmetric form 
\[
\buE(\bv)_t + \sum_{j=1}^d \bffE (\bu(\bv))_{x_j}=\big(A_0^S\big)^{-1}\bv_t + \sum_{j=1}^d A_j^S \bv_{x_j}=0,
\]
with  symmetric Jacobians,  \eqref{eq:symmetry-in-v},  $\displaystyle A_j^S, \ j=0,1,\ldots, d$.
The temporal and spatial fluxes, $\buE(\bv)$ and $\bffE(\bu(\bv))$,  are homogeneous of degree 
$\displaystyle \frac{1}{\beta_\alpha-1} =\frac{\alpha+\gamma}{1-\gamma}<0$.
Since the fluxes are homogeneous of degree $\displaystyle \frac{1}{\beta_\alpha-1}$, Euler's identity enables us to rewrite \eqref{eq:Euler-homo} in the  symmetric conservative form
\[
\big(\big(A_0^S\big)^{-1}\bv\big)_t + \sum_{j=1}^d \big(A_j^S \bv\big)_{x_j}=0,
\]
The homogeneity of Euler equations implies that 
if $\bv(t,\cdot)$ is an entropy variable  solution of Euler equations then any scalar multiple of $\bv$ is also a solution of Euler equations.
What other systems share this property? the system of MHD equations, a natural candidate, does not share this property because of the dependence of its fluxes on dimensional parameter of magnetic permeability  $\mu^*$, \cite[\S5]{Tad1984}. The other system is ultra-relativistic Euler equations which express the conservation of momentum ${\mathbf m}=p \sfv$ and energy $E=3p+p|\sfv|^2$, e.g., \cite{Kun2005,KKMR2024}

\[
\begin{split}
(4m_j\sqrt{1+|\sfv|^2})_t+ \sum_{j=1}^d \big(4m_k \sfv_j +p\delta_{jk}\big)_{x_j}&= 0 \\ E_t + \sum_{j=1}^d \big(4m_k \sqrt{1+|\sfv|^2}\big)_{x_j}&= 0
 \end{split} 
 \]
and we observe that both, the temporal and spatial fluxes are homogeneous of degree 1 in $(p,{\mathbf m})$ (since $\sfv$ is  homogeneous of degree 0).

\medskip\noindent
{\bf Acknowledgment}. 
I benefited from conversations with Francois Golse, Maria Luk\'{a}cov\'{a}-Medvidov\'{a}, Simon Markfleder and Ferdinand Thein. A preliminary version of these results was presented in my SIAM invited address at the 2014 Joint Mathematical Meeting  in Baltimore.


\begin{thebibliography}{1}

\bibitem[AW2021]{AW2021}
I. Akramov and E. Wiedemann,  Non-unique admissible weak solutions of the compressible
Euler equations with compact support in space SIAM J. Math. Anal. 53 (2021) 795-812.

\bibitem[BFH2020]{BFH2020}
D. Breit, E. Feireisl and M. Hofmanová, Solution semiflow to the isentropic Euler system, Arch. Ration. Mech. Anal. 235.1 (2020), pp. 167-194.

\bibitem[Bre2000]{Bre2000} A. Bressan, Hyperbolic systems of conservation laws: the one-dimensional Cauchy problem, OUP Oxford; 2000.

\ifx
\bibitem[CFMP2013]{CFMP2013}
M.J. Castro, U.S. Fjordholm, S. Mishra and C.  Parés, Entropy conservative and entropy stable schemes
for nonconservative hyperbolic systems, SIAM J. Numer. Anal. 51(3), 1371-1391 (2013)
\fi

\bibitem[CL2000]{CL2000}
G. Q. Chen and P. G.  LeFloch, Compressible Euler Equations with general pressure law, Archive Rat. mech. Anal. 153(3) (2000) 221-259.

\bibitem[Che2005]{Che2005}
G.-Q. Chen, Euler equations and related hyperbolic
conservation laws, In ``Handbook of differential equations: evolutionary equations'' 2005, Vol. 2, pp. 1-104,. North-Holland.

\bibitem[CC2007]{CC2007}
G.Q.  Chen and J. Chen, Stability of rarefaction waves and vacuum states for the multidimensional Euler equations. J. Hyperbolic Differ. Equ. 4(1), (2007) 105-122.

\bibitem[Chi2014]{Chi2014}
E. Chiodaroli,  A counterexample to well-posedness of entropy solutions to the compressible Euler system. J. Hyperbolic Differ. Equ. 11(3) (2014) 493-519. 

\bibitem[CDK2015]{CDK2015}
E. Chiodaroli, C. De Lellis, and O. Kreml, Global ill-posedness of the isentropic system of
gas dynamics, Comm. Pure Appl. Math. 68 (2015), pp. 1085-1283.

\bibitem[CK2014]{CK2014}
E. Chiodaroli \& O. Kreml,
On the energy dissipation rate of solutions
to the compressible isentropic Euler system,
Arch. Rational Mech. Anal. 214 (2014) 1019-1049.

\bibitem[CK2018]{CK2018}
E. Chiodaroli and O. Kreml,  Non-uniqueness of admissible weak solutions to the Riemann problem for isentropic Euler equations, Nonlinearity 31 (2018) 1441-60.

\bibitem[CKMS2021]{CKMS2021}
E. Chiodaroli, O. Kreml, V. M\'{a}cha and S.  Schwarzacher, Non-uniqueness of admissible weak solutions to the compressible Euler equations with smooth initial data. Trans. Am. Math. Soc. 374(4), 2269-2295, 2021.

\bibitem[CCT2026]{CCT2026}
J.A. Carrillo, Y.-P. Choi and E. Tadmor,
Lagrangian formulation and Eulerian closure in alignment dynamics,  ArXiv:2604.10253 (2026).

\bibitem[Daf1973]{Daf1973} 
C. Dafermos,
The Entropy Rate Admissibility Criterion for Solutions
of Hyperbolic Conservation Laws, J. Diff. Eq. 14, 202-212 (1973)

\bibitem[Daf1979]{Daf1979}
C. M. Dafermos, The second law of thermodynamics and stability,  Archive for Rational Mechanics and Analysis, 70(2), (1979) 167-179. 

\bibitem[Daf2009]{Daf2009}
C. Dafermos, A variational approach to the Riemann problem for hyperbolic conservation laws, Discrete Cont. Dyn. Systems 23(1\&2), (2009),  185-195. 

\bibitem[Daf2016]{Daf2016}
C. Dafermos, Hyperbolic conservation laws in continuum physics. Vol. 3. Berlin, Springer, 2016.

\bibitem[DLM1995]{DLM1995} G. Dal Maso G, P. LeFloch  and F. Murat, Definition and weak stability of nonconservative products, J. Math. Pures Appl., 74 (1995), pp. 483-548.

\bibitem[DS2010]{DS2010}
C. De Lellis and L. Sz\'{e}kelyhidi, On admissibility criteria for weak solutions of the Euler equations Arch. Ration. Mech. Anal. 195 (2010) 225-60.

\bibitem[DiP1983]{DiP1983}
R. J. DiPerna, Convergence of the viscosity method of isentropic gas dynamics. Commun. Math. Phys. 91, (1983) 1-30.

\bibitem[DiP1985]{DiP1985}
R. J. DiPerna, 
Measure-valued solutions to conservation laws,
Arch. Rat. Mech. Anal. 88 (1985)   223-270.

\bibitem[DLT2026]{DLT2026}
M. Dumbser, M. Luk\'{a}cov\'{a}-Medvidov\'{a} and A. Thomann,
Convergence of a hyperbolic thermodynamically
compatible finite volume scheme for the Euler equations,
\href{https://doi.org/10.1007/s00211-025-01522-2}{Numerische Mathematik  158, (2026) 715-747}.

\bibitem[EL2024]{EL2024}
T. Eiter and R. Lasarzik, Existence of energy-variational solutions to hyperbolic conservation laws, Calc. Var. Partial Differential Eqs  63(4):Paper No. 103, 40, (2024).

\bibitem[FJL2025]{FJL2025}
E. Feireisl, A. J\"{u}ngel, M. Luk\'{a}\v cov\'{a}-Medvid’ov\'{a},
Maximal dissipation and well-posedness of the
Euler system of gas dynamics, ArXiv:2501.05132v2.

\bibitem[FK2015]{FK2015}
E. Feireisl and O. Kreml,  Uniqueness of rarefaction waves in multidimensional compressible Euler system. J. Hyperbolic Differ. Equ. 12(3), (2015) 489-499. 

\bibitem[FL2025]{FL2025}
E. Feireisl, M. Luk\'{a}\v cov\'{a}-Medvid’ov\'{a},
Well-posedness of the Euler system of gas dynamics, ArXiv:2512.18267v1, (2025).

\bibitem[FLY2025]{FLY2025}
E. Feireisl, M. Luk\'{a}\v cov\'{a}-Medvid’ov\'{a},  C. Yu, 
Oscillatory approximations and maximum
entropy principle for the Euler system of gas
dynamics, ArXiv:2505.02070v1.

\ifx%
\bibitem[FMT2012]{FMT2012}
U. S. Fjordholm, S. Mishra and E. Tadmor, Arbitrarily high-order accurate entropy stable essentially
nonoscillatory schemes for systems of conservation laws. SIAM J. Numer. Anal. 50(2), (2012) 544-573.
\fi%

\bibitem[FMT2016]{FMT2016}
U. S. Fjordholm, S. Mishra and E. Tadmor, On the computation of measure-valued solutions, Acta Numer. 25 (2016), pp. 567-679.

\bibitem[Fri1958]{Fri1958} K. O. Friedrichs, Symmetric positive linear differential equations. Comm. Pure Appl. Math. 11 (1958) 333-418.



\bibitem[Fri1980]{Fri1980}
K. O. Friedrichs, Von Neumann’s Hilbert space theory and Partial Differential Equations, SIAM Review 22(4) (1980) 486-493.

\bibitem[FL1971]{FL1971} K. O. Friedrichs  and P. D. Lax, System of conservation equations with a convex extension, Proc. NAS, USA 68(8) 1971,  1686-1688.

\bibitem[Gal2023]{Gal2023}
D. Gallenmüller, Measure-valued low Mach number limits of ideal fluids, SIAM J. Math. Anal. 55.2 (2023), pp. 1145-1169.

\bibitem[GW2021]{GW2021}
D. Gallenmüller and E. Wiedemann, On the selection of measure-valued solutions for the isentropic Euler system,  J. Differential Equations 271 (2021), pp. 979-1006.

\bibitem[GGKS2025]{GGKS2025}
H. Gimperlein, M. Grinfeld, R. J. Knops \& M. Slemrod,
The least action admissibility principle,
Arch. Rational Mech. Anal. (2025) 249:22.

\bibitem[GK2022]{GK2022}
 V. Giri and  H. Kwon, On non-uniqueness of continuous entropy solutions to the Isentropic Compressible Euler Equations, Arch. Rational Mech. Anal. 245 (2022) 1213-1283.
 
\bibitem[Gli1965]{Gli1965} J. Glimm, 
Solutions in the large for nonlinear hyperbolic systems of equations (1965)
Communications on Pure and Applied Math. 18(4) (1965) 697-715.
  
\bibitem[GLC2020]{GLC2020}
 J. Glimm, D.  Lazarev and G.-Q. Chen, Maximum entropy production as a necessary admissibility condition for the fluid Navier–Stokes and Euler equations. SN Appl. Sci. 2(2160), 2020.
 
\bibitem[God1961]{God1961} S. K. Godunov, An interesting class of quasilinear systems, Dokl. Akad. Nauk  139(3) (1961), 521-523 (translated in \href{https://doi.org/10.1016/j.jcp.2024.113521}{Journal of Computational Physics 520 (2025) 113521}).
 
  \bibitem[God1962]{God1962}
 S. K. Godunov, The problem of a generalized solution in the theory of quasilinear equations and in gas dynamics,  Russian Mathl Surveys 17(3) (1962) 145.
 
\bibitem[God1972]{God1972}
S. K. Godunov, Symmetric form of the equations of magnetohydrodynamics. Numerical Methods for
Mechanics of Continuous Media 3(1) (1972) 26-31.

\bibitem[God1986]{God1986}
S. K. Godunov, Lois de conservation et int\'{e}grales d'energie des equations hyperboliques,
in ``Nonlinear Hyperbolic Problems'', Proceedings of a 1986 Advanced Research Workshop, Lecture Notes in Mathematics, vol. 1270 (C. Carasso, P.-A. Raviart and D. Serre, eds.), Springer-Verlag, 1987, pp-135--149.

\bibitem[Har1983]{Har1983}
A. Harten, On the symmetric form of systems of conservation laws with entropy, \href{https://doi.org/10.1016/0021-9991(83)90118-3}{J. Comput. Phys. 49(1), (1983) 151-164}. 

\bibitem[HHW2020]{HHW2020}
Y. Holle, M. Herty and M. Westdickenberg,
New coupling conditions for isentropic flow on networks, Networks and Homogeneous Media, 15(4) (2020), 605-631

\bibitem[Hsi1980]{Hsi1980}
L. Hsiao, The entropy rate admissibility criterion for gas dynamics, J. Differential
Equations 38 (1980), 226-238.

\bibitem[HFM1986]{HFM1986}
T.J.R. Hughes, L.P. Franca and M. Mallet, Symmetric forms of the compressible Euler and Navier-Stokes equations and the second law of thermodynamics, Comput. Methods. Appl. Mech. Engrg. 54 (1986) 223-234. 

\bibitem[Jen2000]{Jen2000} H. K. Jenssen, 
Blowup for Systems of Conservation Laws,
SIAM J. on Math. Anal. 31(4) (2000) 894-908.

\bibitem[Kle2023]{Kle2023}
S.-C. Klein,
Stabilizing discontinuous Galerkin methods using Dafermos' entropy rate criterion: I-one-dimensional conservation laws,
J. Scientific Computing 95(55) (2023) 1-37.

\bibitem[KKMM2020]{KKMM2020}
C. Klingenberg, O. Kreml, V. M\'{a}cha, and S. Markfelder, Shocks make the Riemann problem
for the full Euler system in multiple space dimensions ill-posed, Nonlinearity 33.12 (2020),  6517-6540.

\bibitem[KMW2026]{KMW2026}
C. Klingenberg, S. Markfelder and  E. Wiedemann,
Maximal turbulence as a selection criterion for measure-valued solutions
ArXiv:2503.20343 (2026).

\bibitem[KV2019]{KV2019}
S. G. Krupa and A. F. Vasseur,
On uniqueness of solutions to conservation laws verifying a single entropy condition, J. Hyperbolic Differential Eq. 16(1) (2019) 157-191.


\bibitem[Kru1970]{Kru1970}
S. Kru\v{z}kov, First order quasilinear equations in several independent variables, Math. of the USSR-Sbornik 10.2 (1970): 217.

\bibitem[Kun2005]{Kun2005}
M. Kunik, 
Selected Initial and Boundary Value Problems for Hyperbolic Systems and Kinetic Equations, Habilitation thesis, Otto-von-Guericke University Magdeburg, 2005, The thesis is available under \href{https://opendata.uni-halle.de//handle/1981185920/30710}{https://opendata.uni-halle.de//handle/1981185920/30710}.

\bibitem[KKMR2024]{KKMR2024}
M. Kunik, A. Kolb, S. M\"uller and H. Ranocha, Radially symmetric solutions of the ultra-relativistic Euler equations in several space dimensions, J. Comput. Phys., 518 (2024), 113330.

\bibitem[Lax1957]{Lax1957}
P. D. Lax, Hyperbolic systems of conservation laws. Comm. Pure Appl. Math. 10 (1957), 537-566.

\bibitem[Lax1971]{Lax1971}
P.  Lax, Shock waves and entropy, in ``\emph{Contributions to nonlinear functional analysis}'', Academic Press, 1971. 603-634.


\bibitem[Lax1973]{Lax1973}
P. Lax,  Hyperbolic systems of conservation laws and the mathematical theory of shock waves, {SIAM}, 1973.

\bibitem[Lax1987]{Lax1987} P. Lax, On symmetrizing hyperbolic differential equations, in ``Nonlinear Hyperbolic Problems'', Proceedings of a 1986 Advanced Research Workshop, Lecture Notes in Mathematics, vol. 1270 (C. Carasso, P.-A. Raviart and D. Serre, eds.), Springer-Verlag, 1987, pp. 150--151.

\bibitem[Lax1989]{Lax1989}
P. D. Lax, 
The Flowering of Applied Mathematics in America,
SIAM Review, 31(4) (1989) 533-541.

\bibitem[Lax2005]{Lax2005}
P. D. Lax,
Selected Papers I and II (Peter Sarnak and Andrew J. Majda, eds),
Springer 2005.

\bibitem[Lax2007]{Lax2007}
P. D. Lax, Mathematics and Physics. Bull. AMS 45(1), 2007, 135-152.

\bibitem[Lax2014]{Lax2014}
P. D. Lax, 
\href{https://www.math.umd.edu/\~tadmor/ki\_net/activities/tn60/2014\_04\_30\_Lax\_Banquet\_talk.pdf}{John von Neumann: The Early Years, The Years at Los Alamos
and the Road to Computing}, in ``\href{https://www.math.umd.edu/\~tadmor/ki\_net/activities/tn60/}{\textcolor{black}{Modern Perspectives in Applied Mathematics: Theory and Numerics of PDEs}}'', College Park, 2014.


\bibitem[LOW2004]{LOW2004}
C. De Lellis, F. Otto, and M. Westdickenberg, Minimal entropy conditions for Burgers equation. Quart. Appl. Math., 62(4) (2004) 687-700.

\bibitem[LS2010]{LS2010}
C. De Lellis and L. Székelyhidi, On admissibility criteria for weak solutions of the Euler equations. Arch. Ration. Mech. Anal. 195(1), (2010) 225-260.

\bibitem[LPT1994]{LPT1994}
P.-L. Lions, P. Perthame and E. Tadmor, 
Kinetic formulation of the isentropic gas dynamics and p-systems
Communications in Mathematical Physics 163 (1994), 415-431.

\bibitem[Mar2024]{Mar2024}
S. Markfelder, A new convex integration approach for the
compressible Euler equations and failure of
the local maximal dissipation criterion,
Nonlinearity 37 (2024) 115022 1-60.

\bibitem[MK2018]{MK2018}
S. Markfelder and C. Klingenberg,
The Riemann problem for the
multidimensional isentropic system of gas
dynamics is ill-posed if it contains a shock,
Arch. Rational Mech. Anal. 227 (2018) 967-994.

\bibitem[MS2006]{MS2006}
L. M. Martyushev and V. D.  Seleznev, Maximum entropy production principle in physics, chemistry and biology, Physics Reports 426 (2006) 1- 45.

\bibitem[Moc1980]{Moc1980}
M. S. Mock, Systems of conservation laws of mixed type,  Journal of Differential Equations. 37(1) (1980), 70-88.

\bibitem[Pan1994]{Pan1994} E. Panov, Uniqueness of the solution of the Cauchy 
problem for a first-order quasilinear equation with an admissible strictly convex entropy,  Mat. 
Zametki  55(5)  (1994),  116-129, 159 (translation in  Math. Notes  55  (1994),  no. 5-6, 517-525).

 \bibitem[Roz1959]{Roz1959}
 B. L. Rozhdestvenskii, On the conservativeness of systems of quasilinear
equations, Uspekhi Mat. Nauk 14 (1959), 217-218. (Russian)

\bibitem[Roz1960]{Roz1960}
B. L. Rozhdestvenskii, Discontinuities solutions of hyperbolic systems of quasilinear equations,  Russ. Math. Surv. 15 (1960) 53-111.
 
 \bibitem[Ser1999]{Ser1999}
 D. Serre,  Systems of Conservation Laws, Vol 2. Geometric Structures, Oscillations, and Initial-Boundary Value Problems. Cambridge University Press; 1999.
 
\bibitem[Sev1990]{Sev1990}
M. Sever, The rate of total entropy generation for Riemann problems, J. Differential Equations 87 (1990), 115-143.

\bibitem[Sev1996]{Sev1996}
M. Sever, Estimate of the time rate of entropy dissipation
for systems of conservation laws,
J. Differential Equations 130, (1996) 127-141.

\bibitem[Tad1984]{Tad1984}
E. Tadmor, 
Skew selfadjoint form for systems of conservation laws,
Journal of Mathematical Analysis and Applications 103(2) (1984), 428-442

\ifx
\bibitem[Tad1986]{Tad1986}
E. Tadmor, 
A minimum entropy principle in the gas dynamics equations
Applied Numerical Mathematics 2 (1986), 211-219.
\fi

\bibitem[Tad1987a]{Tad1987a}
E. Tadmor, 
The entropy dissipation by numerical viscosity in nonlinear conservative difference schemes,
in ``Nonlinear Hyperbolic Problems'', Proceedings of a 1986 Advanced Research Workshop, Lecture Notes in Mathematics, vol. 1270 (C. Carasso, P.-A. Raviart and D. Serre, eds.), Springer-Verlag, 1987, pp. 52--63.


\bibitem[Tad1987b]{Tad1987b}
E. Tadmor, The numerical viscosity of entropy stable schemes for systems of conservation laws I.
Math. Comput. 49 (1987) 91-103.

\bibitem[Tad1987c]{Tad1987c}
E. Tadmor, 
Entropy functions for symmetric systems of conservation laws,
Journal of Mathematical Analysis and Applications 122(2) (1987) 355-359.

\bibitem[Tad2003]{Tad2003}
E. Tadmor, Entropy stability theory for difference approximations of nonlinear conservation laws and related time dependent problems,
\href{https://doi.org/10.1017/S0962492902000156}{Acta Numerica 12 (2003) 451-512}.

\ifx
\bibitem[Tad2021]{Tad2021}
E. Tadmor, On the mathematics of swarming: emergent behavior in alignment dynamics, 
Notices of the AMS 68(4) (2021) 493-503.

\bibitem[Tad2023]{Tad2023} 
E. Tadmor, Swarming: hydrodynamic alignment with pressure, Bulletin AMS 60(3) (2023) 285-325.
\fi

\bibitem[Tad2025]{Tad2025}
E. Tadmor, Entropy decrease and emergence of order in collective dynamics, Communications in Contemporary Math. 28(5) (2025) 2540006 1-23.

\ifx
\bibitem[TR2025]{TR2025}
F.  Thein  and H. Ranocha,
Computing radially-symmetric solutions of the ultra-relativistic Euler equations with entropy-stable discontinuous Galerkin methods,
ArXiv:2508.21427 (2025).
 \fi
 
\bibitem[Vol1967]{Vol1967} A. I. Volpert,   The spaces BV and quasilinear equations. Mat. Sbornik, 73(115), (1967), No. 2 (1967) 225-267.

\bibitem[wiki:Martians]{wiki:Martians}
Wikipedia contributors, ``The Martians (scientists),'' \textit{Wikipedia, The Free Encyclopedia}, 
\url{https://en.wikipedia.org/wiki/The_Martians_(scientists)} 
(accessed May 21, 2026).
\end{thebibliography}
\end{document}